\documentclass[final]{siamltex}
\usepackage{graphicx,amssymb,psfrag}

\newcommand{\BEAS}{\begin{eqnarray*}}
\newcommand{\EEAS}{\end{eqnarray*}}
\newcommand{\BEA}{\begin{eqnarray}}
\newcommand{\EEA}{\end{eqnarray}}
\newcommand{\BEQ}{\begin{equation}}
\newcommand{\EEQ}{\end{equation}}
\newcommand{\BIT}{\begin{itemize}}
\newcommand{\EIT}{\end{itemize}}
\newcommand{\BNUM}{\begin{enumerate}}
\newcommand{\ENUM}{\end{enumerate}}

\newcommand{\BA}{\begin{array}}
\newcommand{\EA}{\end{array}}

\newcommand{\ones}{\mathbf 1}

\newcommand{\reals}{{\mbox{\bf R}}}

\newcommand{\symm}{{\mbox{\bf S}}}  

\newcommand{\Card}{\mathop{\bf Card}}
\newcommand{\Tr}{\mathop{\bf Tr}}

\newcommand{\id}{{\mathbf I}}

\newcommand{\argmin}{\mathop{\rm argmin}}

\title{First-order methods for sparse covariance selection}

\author{
Alexandre d'Aspremont\thanks{ORFE Department, Princeton University, Princeton, NJ 08544. \texttt{\small{aspremon@princeton.edu}}}
\and
Onureena Banerjee\thanks{EECS Department, UC Berkeley, Berkeley, CA 94720. \texttt{\small{onureena@eecs.berkeley.edu}}}
\and
Laurent El Ghaoui\thanks{EECS Department, UC Berkeley, Berkeley, CA 94720. \texttt{\small{elghaoui@eecs.berkeley.edu}}}}
\begin{document}

\maketitle

\begin{abstract}
Given a sample covariance matrix, we solve a maximum likelihood problem penalized by the number of nonzero coefficients in the inverse covariance matrix. Our objective is to find a sparse representation of the sample data and to highlight conditional independence relationships between the sample variables. We first formulate a convex relaxation of this combinatorial problem, we then detail two efficient first-order algorithms with low memory requirements to solve large-scale, dense problem instances.
\end{abstract}

\section{Introduction}
We discuss a problem of model selection\footnote{A subset of the results discussed here appeared in the proceedings of the International Conference on Machine Learning, Pittsburgh 2006.}. Given $n$ variables drawn from a Gaussian distribution $\mathcal{N}(0,C)$, where the true covariance matrix $C$ is unknown, we estimate $C$ from a sample covariance matrix $\Sigma$ by maximizing its log-likelihood. Following \cite{Demp72}, setting a certain number of coefficients in the inverse covariance matrix $\Sigma^{-1}$ to zero, a procedure known as \emph{covariance selection}, improves the stability of this estimation procedure by reducing the number of parameters to estimate and highlights structure in the underlying model.

Here, we focus on the problem of \emph{discovering} this pattern of zeroes in the inverse covariance matrix. We seek to trade-off the log-likelihood of the solution with the number of zeroes in its inverse, and solve the following estimation problem:
\BEQ \label{eq:sparseml-primal}
\BA{ll}
\mbox{maximize} & \log \det X - \langle \Sigma, X \rangle - \rho \Card(X)\\
\mbox{subject to} & \alpha \id_n \preceq X \preceq \beta \id_n
\EA
\EEQ
in the variable $X\in \symm_n$, where $\Sigma \in \symm_n^+$ is the sample covariance matrix, $\Card(X)$ is the cardinality of $X$, i.e. the number of nonzero components in $X$, $\rho >0$ is a parameter controlling the trade-off between log-likelihood and cardinality, finally $\alpha,~\beta >0$ fix bounds on the eigenvalues of the solution.

Zeroes in the inverse covariance matrix correspond to conditionally independent variables in the model and this approach can be used to simultaneously determine a robust estimate of the covariance matrix and, perhaps more importantly, discover \emph{structure} in the underlying graphical model. In particular, we can view (\ref{eq:sparseml-primal}) as a model selection problem using Aikake (AIC, see \cite{Akai73}) or Bayes (BIC, see \cite{Burn04}) information criterions. Both these problems can be written as in (\ref{eq:sparseml-primal}) with $\rho = 2/N$ for the AIC problem and $\rho=2\log(N/2)/N$ for the BIC problem, where $N$ is the sample size. This has applications in speech recognition (see \cite{bilmes1999,bilmes2000}) or gene networks analysis (see \cite{Dobr04,Dobr04a} for example).


The $\Card(X)$ penalty term makes the estimation problem (\ref{eq:sparseml-primal}) combinatorial (NP-Hard in fact), and our first objective here is to derive a convex relaxation to this problem which can be solved efficiently. We then derive two first-order algorithms geared towards memory efficiency and large-scale, \emph{dense} problem instances.

In \cite{bilmes2000}, Bilmes proposed a method for covariance selection based on choosing statistical dependencies according to conditional mutual information computed using training data. Other recent work involves identifying those Gaussian graphical models that are best supported by the data and any available prior information on the covariance matrix. This approach is used by \cite{jones2004,Dobr04} on gene expression data. Recently, \cite{Dahl05,Huan05} also considered penalized maximum likelihood estimation for covariance selection. In contrast to our results here, \cite{Huan05} work on the Cholesky decomposition of $X$ using an iterative (heuristic) algorithm to minimize a nonconvex penalized likelihood problem, while \cite{Dahl05} propose a set of large scale interior point algorithms to solve sparse problems, i.e. problems for which the conditional independence structure is already known.

The paper is organized as follows, in Section \ref{s:setup}, we detail our convex relaxation of problem (\ref{eq:sparseml-primal}) and study the dual. In Section \ref{s:algos}, we derive two efficient algorithms to solve it. Finally, in Section \ref{s:numer} we describe some numerical results.

\section{Problem setup}\label{s:setup}
\subsection{Convex relaxation}
Given a sample covariance matrix $\Sigma\in\symm^+_n$, we can write the following convex relaxation to the estimation problem (\ref{eq:sparseml-primal}):
\BEQ \label{eq:sparseml-relax}
\BA{ll}
\mbox{maximize} & \log \det X - \langle \Sigma, X \rangle - \rho \ones^T|X|\ones\\
\mbox{subject to} & \alpha \id_n \preceq X \preceq \beta \id_n,
\EA
\EEQ
with variable $X \in \symm^n$, where $\ones$ is the $n$-vector of ones, so that $\ones^T|X|\ones =  \sum_{i,j=1}^n |X_{ij}|$.  
The penalty term involving the sum of absolute values of the entries of $X$ is a proxy for the number of its non-zero elements: the function $\ones^T|X|\ones$ can be seen as the largest convex lower bound on $\Card(X)$ on the hypercube, an argument used by \cite{Boyd00} for rank minimization.
It is also often used in regression techniques, such as the LASSO studied by \cite{Tibs96}, when sparsity of the solution is a concern. This relaxation is provably tight in certain cases (see \cite{Dono05}). In our model, the bounds $(\alpha,\beta)$ on the eigenvalues of $X$ are fixed, and user-chosen. Although we allow $\alpha = 0$, $\beta = +\infty$, such bounds are useful in practice to control the condition number of the solution.

For $\rho = 0$, and provided $\Sigma \succ 0$, problems (\ref{eq:sparseml-primal}) and (\ref{eq:sparseml-relax}) have a unique solution $X^\star = \Sigma^{-1}$, and the corresponding maximum-likelihood estimate is $\Sigma$. Due to noise in the data, in practice, the sample estimate $\Sigma$ may not have a sparse inverse, even if the underlying graphical model exhibits conditional independence properties. By striking a trade-off between maximality of the likelihood and number of non-zero elements in the inverse covariance matrix, our approach is potentially useful at \emph{discovering structure}, precisely conditional independence properties in the data. This means that we have to focus on the case where the matrix $X$ is \emph{dense}. At the same time, it serves as a regularization technique: when $\Sigma$ is rank-deficient, there is no well-defined maximum-likelihood estimate, whereas the solution to problem
(\ref{eq:sparseml-relax}) is always unique and well-defined for $\rho>0$, as seen below.

\subsection{Dual problem, robustness}
We can rewrite the relaxation (\ref{eq:sparseml-relax}) as the following min-max problem:
\BEQ \label{eq:sparseml-robust}
\max_{\{X:~\alpha \id_n \preceq X \preceq \beta \id_n\}} \min_{\{U:~|U_{ij}|\leq \rho\}} \log \det X - \langle \Sigma + U, X \rangle\\
\EEQ
which gives a natural interpretation of problem (\ref{eq:sparseml-relax}) as a worst-case \emph{robust maximum likelihood} problem with componentwise bounded, additive noise on the sample covariance matrix $\Sigma$. The corresponding Lagrangian is given by:
\[
L(X,U,P,Q)=\log\det X -\Tr((\Sigma+U+Q-P)X) - \alpha \Tr P + \beta \Tr Q
\]
and we get the following dual to (\ref{eq:sparseml-relax}):
\BEQ \label{eq:sparseml-dual}
\BA{ll}
\mbox{minimize} & -\log \det (\Sigma + U + Q - P) -n + \alpha \Tr P - \beta \Tr Q\\
\mbox{subject to} & P,Q \succeq0,~|U_{ij}|\leq \rho, \quad i,j=1,\ldots,n,
\EA
\EEQ
in the variables $U,P,Q\in\symm_n$. 

When $\alpha=0$ and $\beta=+\infty$, the first-order optimality conditions impose $X(\Sigma + U)=\id_n$, hence we always have:
\[
X \succeq \alpha(n) \id_n \quad \mbox{with} \quad \alpha(n):=\frac{1}{\|\Sigma\|+n\rho},
\]
zero duality gap also means $\Tr (\Sigma X)= n- \rho \ones^T|X|\ones$. Because $X$ and $\Sigma$ are both positive semidefinite, we get:
\[
\|X\|_2 \leq \|X\|_F \leq \ones^T|X|\ones \leq \frac{n}{\rho},
\]
which, together with $\Tr(\Sigma X) \geq \lambda_\mathrm{min}(\Sigma) \|X\|_2$, means $\|X\|_2 \leq n/\lambda_\mathrm{min}(\Sigma)$. Finally then, we must always have:
\[
X \preceq \beta(n) \id_n \quad \mbox{with} \quad \beta(n):=n \min\left(\frac{1}{\rho},\|\Sigma^{-1}\|_2\right).
\]
and $ 0< \alpha(n) \leq \lambda(X) \leq \beta(n) < + \infty$ at the optimum. Setting $\alpha=0$ and $\beta=+\infty$ in problem (\ref{eq:sparseml-relax}) is then equivalent to setting $\alpha=\alpha(n)$ and $\beta=\beta(n)$. Since the objective function of problem (\ref{eq:sparseml-relax}) is strictly convex when $0< \alpha(n) \leq \lambda(X) \leq \beta(n) < + \infty$, this shows that (\ref{eq:sparseml-relax}) always has a unique solution.

\section{Algorithms} \label{s:algos}
In this section, we present two algorithms for solving problem (\ref{eq:sparseml-relax}), one based on an optimal first-order method developed in \cite{Nest03}, the other based on a block-coordinate gradient method.

Of course, problem (\ref{eq:sparseml-relax}) is convex and can readily be solved using interior point methods (see \cite{Boyd03} for example). However, such second-order methods become quickly impractical for solving (\ref{eq:sparseml-primal}), since the corresponding complexity to compute an $\epsilon$-suboptimal solution is $O(n^6 \log (1/\epsilon))$. Note however that we cannot expect to do better than $O(n^3)$, which is the cost of solving the non-penalized problem for dense covariance matrices $\Sigma$. 
\subsection{Smooth optimization} \label{ss:nesterov}
The recently-developed first-order algorithms due to \cite{Nest03} trade-off a better dependence on problem size against a worst dependence on accuracy, usually $1/\epsilon$ instead of its logarithm and the method we describe next has a complexity of $O(n^{4.5}/\epsilon)$. In addition, the memory space requirement of these first-order methods is much lower than that of interior-point methods, which involve forming a dense Hessian, and hence become quickly prohibitive with a problem having $O(n^2)$ variables. 

\paragraph{Nesterov's format} 
The algorithm in \cite{Nest03} supposes that the function to minimize conforms to a certain representation. This is the case for our problem here, so we first write (\ref{eq:sparseml-robust}) in the saddle-function format described in \cite{Nest03}:
\[
\min_{X \in {\cal Q}_1} \: -\log \det X + \langle \Sigma, X \rangle + \rho \ones^T|X|\ones \equiv \min_{X \in {\cal Q}_1} \: \max_{U \in {\cal Q}_2} \: \hat{f}(X) + \langle A(X),U \rangle
\]
where we define $\hat{f}(X) = -\log \det X +  \langle \Sigma,X \rangle$, $A = \rho I_{n^2}$, and
\[
{\cal Q}_1 := \left\{ X \in {\cal S}^n ~:~ \alpha \id_n \preceq X \preceq \beta \id_n \right\}, \quad {\cal Q}_2 := \left\{ U \in {\cal S}^n  ~:~ \|U\|_\infty \leq 1 \right\}.
\]
The adjoint of this problem, corresponding to the dual problem (\ref{eq:sparseml-dual}), is then written:
\BEQ \label{eq:sparseml-adjoint}
\max_{U \in {\cal Q}_2} \phi(U) \quad \mbox{where} \quad \phi(U):= \min_{X \in {\cal Q}_1} - \log \det X + \langle \Sigma + U, X \rangle.
\EEQ
When a function can be represented in this saddle function format, the method described in \cite{Nest03} combines two steps. \emph{Regularization}: by adding a strongly convex penalty to the saddle function representation of $f$, the algorithm first computes a smooth $\epsilon$-approximation of $f$ with Lipschitz continuous gradient. This can be seen as a generalized Moreau-Yosida regularization step (see \cite{Lema97} for example). \emph{Optimal first order minimization}: the algorithm then applies the optimal first-order scheme for functions with Lipschitz continuous gradients detailed in \cite{Nest83} to the regularized function. Each iteration requires efficiently computing the regularized function value and its gradient. In all the semidefinite programming applications detailed here, this can be done extremely efficiently, with a complexity of $O(n^3)$ and memory requirements in $O(n^2)$. The method is only efficient if all these steps can be performed explicitly or at least very efficiently. As we will see below, this is the case here.

\paragraph{Prox-functions and related parameters}
To ${\cal Q}_1$ and ${\cal Q}_2$ we now associate norms and so-called prox-functions. For ${\cal Q}_1$, we use the Frobenius norm, and a prox-function:
\[
d_1(X) = -\log \det X + \log \beta .
\]
The function $d_1$ is strongly convex on ${\cal Q}_1$, with a convexity parameter of $\sigma_1 =1/\beta^2$, in the sense that $\nabla^2 d_1(X)[H,H] = \Tr (X^{-1}H X^{-1}H) \ge \beta^{-2} \|H\|_F^2$ for every $H$. Furthermore, the center of the set, $X_0:=\arg\min_{X \in {\cal Q}_1} \: d_1(X)$ is $X_0 = \beta \id_n$, and satisfies $d_1(X_0) = 0$. With our choice, we have $D_1 :=\max_{X \in {\cal Q}_1} \: d_1(X) = n\log (\beta/\alpha)$.

To ${\cal Q}_2$, we also associate the Frobenius norm, and the prox-function $d_2(U) = \|U\|_F^2/2$. With this choice, the center $U_0$ of ${\cal Q}_2$ is $U_0 = 0$. Furthermore, the function $d_2$ is strongly convex on its domain, with convexity parameter with respect to the $1$-norm $\sigma_1 = 1$, and we have $D_2 := \max_{U \in {\cal Q}_2} \: d_2(U) = n^2/2$.

The function $\hat{f}$ has a gradient that is Lipschitz-continuous with respect to the Frobenius norm on the set ${\cal Q}_1$, with Lipschitz constant $M = 1/\alpha^2$. Finally, the norm (induced by the Frobenius norm) of the operator $A$ is $\|A\| = \rho$.

\paragraph{Smooth minimization}
The method is based on replacing the objective of the original problem, $f(X)$, with $f_\epsilon(X)$, where $\epsilon>0$ is the desired accuracy, and $f_\epsilon$ is a penalized function involving the prox-function $d_2$, defined as: 
\BEQ\label{eq:feps-def}
f_\epsilon(X) := \hat{f}(X) + \max_{U \in {\cal Q}_2} \langle X,U \rangle - (\epsilon/2D_2) d_2(U).
\EEQ
The above function turns out to be a smooth uniform approximation to $f$ everywhere, with maximal error $\epsilon/2$.  Furthermore, the function $f_\epsilon$ is has a Lipschitz-continuous gradient, with Lipschitz constant given by $L(\epsilon) := M+D_2\|A\|^2/(2 \sigma_2\epsilon)$. A specific first-order algorithm detailed in \cite{Nest83} for smooth, constrained convex minimization is then applied to the function $f_\epsilon$, with convergence rate in $O(\sqrt{L(\epsilon)/\epsilon})$.

\paragraph{Nesterov's algorithm} Choose $\epsilon>0$ and set $X_0 =  \beta\id_n$, the algorithm then updates primal and dual iterates $Y_k$ and $\hat U_k$ using the following steps:
\vskip 1ex
\begin{enumerate} \itemsep 0ex
    \item Compute $\nabla f_\epsilon(X_k) = -X^{-1} + \Sigma + U^\ast(X_k)$, where $U^\ast(X)$
    solves (\ref{eq:feps-def}).
    \item Find $Y_k = \argmin_{Y  \in {\cal Q}_1} \: \{\langle \nabla f_\epsilon(X_k) , Y-X_k \rangle +
\frac{1}{2} L(\epsilon) \|Y-X_k\|_F^2\}$
    \item Find $Z_k = \argmin_{Z \in {\cal Q}_1} \left\{ \frac{L(\epsilon)d_1(Z)}{\sigma_1}  + \sum_{i=0}^k \frac{i+1}{2}
(f_\epsilon(X_i)+    \langle \nabla f_\epsilon(X_i),Z-X_i \rangle) \right\}$
    \item Update $X_k = \frac{2}{k+3} Z_k + \frac{k+1}{k+3} Y_k$ and $\hat U_k=\frac{k\hat U_{k-1} + 2U^\ast(X_k)}{(k+2)}$
    \item Repeat until the duality gap is less than the target precision:
    \[
    -\log \det Y_k + \langle \Sigma, Y_k \rangle + \rho \ones^T|Y_k|\ones - \phi(\hat U_k) \leq \epsilon.
    \]
\end{enumerate}

\noindent The key to the method's success is that step 1-3 and 5 can be performed explicitly and only involve an eigenvalue decomposition. Step one above computes the (smooth) function value and gradient. The second step computes the \emph{gradient mapping}, which matches the gradient step for unconstrained problems (see \cite[p.86]{Nest03a}). Step three and four update an \emph{estimate sequence} \cite[p.72]{Nest03a} of $f_\mu$ whose minimum can be computed explicitly and gives an increasingly tight upper bound on the minimum of $f_\mu$. We now present these steps in detail for our problem. 

\paragraph{Step 1} The first step requires computing the gradient of the function
\[
f_\epsilon(X) = \hat{f}(X) + \max_{u \in {\cal Q}_2} \: \langle X,U \rangle - (\epsilon/2D_2) d_2(U) .
\]
This function can be expressed in closed form as $f_\epsilon(X) = \hat{f}(X) + \sum_{i,j} \psi_\mu(X_{ij})$, where
\[
\psi_\epsilon (x) := \left\{ \begin{array}{ll} |x| - (\epsilon/4D_2) & \mbox{if } |x| \ge (\epsilon/2\rho D_2), \\
{D_2 x^2}/{\epsilon} & \mbox{otherwise,} \end{array} \right.
\]
which is simply the Moreau-Yosida regularization of the absolute value and the gradient of the function at $X$ is
\[
\nabla f_\mu(X) = -X^{-1} + \Sigma + U^\ast(X) ,
\]
with
\[
U^\ast(X) := \max(\min( X/\mu,\rho),-\rho),
\]
with min. and max. understood componentwise. The cost of this step is dominated by that of computing the inverse of $X$, which is $O(n^3)$.

\paragraph{Step 2}
This step involves a problem of the form
\[
T_{{\cal Q}_1}(X) = \arg\min_{Y \in {\cal Q}_1} \: \langle \nabla f_\epsilon(X) , Y-X \rangle + \frac{1}{2} L \|Y-X\|_F^2 ,
\]
where $X \in {\cal Q}_1$ is given. This problem can be reduced to one of projection on ${\cal Q}_1$, namely
\[
\min_{Y \in {\cal Q}_1} \: \|Y-G\|_F^2 ,
\]
where $G := X - L^{-1}\nabla f_\epsilon(X)$. Using the rotational invariance of this problem, we reduce it to a vector problem:
\[
\textstyle \min_\lambda \: \sum_i (\lambda_i - \gamma_i)^2 ~:~ \alpha \le \lambda_i \le \beta, \quad
i=1,\ldots,n,
\]
where $\gamma$ is the vector of eigenvalues of $G$.  This problem admits a simple explicit solution:
\[
\lambda_i = \min ( \max (\gamma_i,\alpha),\beta ) , \quad i=1,\ldots,n.
\]
The corresponding solution is then $Y = V^T \mbox{\bf diag}(\lambda)V$, where $G = V^T\mbox{\bf diag}(\gamma)V$ is the eigenvalue decomposition of $G$. The cost of this step is dominated by the cost of forming the eigenvalue decomposition of $G$, which is $O(n^3)$.

\paragraph{Step 3} The third step involves solving a problem of the form
\begin{equation}\label{eq:z-def}
Z : = \arg\max_{X \in {\cal Q}_1} \: d_1(X) + \langle S,X \rangle,
\end{equation}
where $S$ is given.  Again, due to the rotational invariance of the objective and feasible set, we can reduce the problem to a one-dimensional problem:
\[
\textstyle \min_\lambda \: \sum_i \sigma_i \lambda_i - \log \lambda_i  ~:~ \alpha \le \lambda_i \le \beta,
\]
where $\sigma$ contains the eigenvalues of $S$. This problem has a simple, explicit solution:
\[
\lambda_i = \min ( \max (1/\sigma_i,\alpha),\beta ) , \quad i=1,\ldots,n.
\]
The corresponding solution is then $Y = V^T \mbox{\bf diag}(\lambda)V$, where $S = V^T\mbox{\bf diag}(\sigma)V$ is the eigenvalue decomposition of $S$. Again, the cost of this step is dominated by the cost of forming the eigenvalue decomposition of $S$, which is $O(n^3)$.

\paragraph{Computing $\phi(\hat U_k)$}
For a given matrix $\hat U_k$ the function $\phi$ is computed as in (\ref{eq:sparseml-adjoint}):
\[
\phi(\hat U_k)= \min_{X \in {\cal Q}_1} - \log \det X + \langle \Sigma + \hat U_k, X \rangle.
\]
This means projecting $(\Sigma + \hat U_k)^{-1}$ on ${\cal Q}_1$, i.e. only involves an eigenvalue decomposition.

\paragraph{Complexity estimate}
To summarize, for step $1$, the gradient of $f_\epsilon$ is readily computed in closed form, via the computation of the inverse of $X$. Step $2$  essentially amounts to projecting on ${\cal Q}_1$, and requires an eigenvalue problem to be solved; likewise for step $3$.  In fact, each iteration costs $O(n^3)$. The number of iterations necessary to achieve an objective with absolute accuracy less than $\epsilon$ is then given by:
\begin{equation}\label{eq:N-def}
N(\epsilon) := 4 \|A\| \frac{1}{\epsilon} \sqrt{\displaystyle\frac{D_1D_2}{\sigma_1\sigma_2} } 
 + \sqrt{\frac{MD_1}{\sigma_1\epsilon}} = \frac{\kappa \sqrt{n(\log
\kappa)}}{\epsilon} ( 4 n \alpha \rho + \sqrt{\epsilon}),
\end{equation}
where $\kappa=\log(\beta/\alpha)$ bounds the solution's condition number. Thus, the overall complexity when $\rho>0$ is in $O(n^{4.5}/\epsilon)$, as claimed.

\subsection{Block-coordinate gradient methods} \label{ss:block}
In this section, we focus on the particular case where $\alpha=0$ and $\beta=+\infty$ (hence implicitly $\alpha=\alpha(n)$, $\beta=\beta(n)$) and derive gradient minimization algorithms that take advantage of the problem structure. We consider the following problem:
\BEQ \label{eq:sparseml-block}
\max_X ~ \log \det X - \langle \Sigma, X \rangle - \rho \ones^T |X| \ones \\
\EEQ 
in the variable $X\in\symm_n$, where $\rho >0$ again controls the trade-off between log-likelihood and sparsity of the inverse covariance matrix. Its dual is given by:
\BEQ \label{eq:dual-sparseml-block}
\BA{ll}
\mbox{minimize} & -\log \det (\Sigma + U )-n\\
\mbox{subject to} & |U_{ij}|\leq \rho, \quad i,j=1,\ldots,n.
\EA\EEQ
in the variable $U\in\symm_n$. We partition the matrices $X$ and $U$ in block format:
\[
X=\left(\BA{cc}
Z & x\\
x^T & y\\
\EA\right)
\quad \mbox{and} \quad
U=\left(\BA{cc}
V & u\\
u^T & w\\
\EA\right),
\]
where $Z \succ 0$ and $U$ are fixed and $x,~u \in \reals^{(n-1)},~ y,~w\in\reals$ are the variables (row and column) we are updating. We also partition the sample matrix according to the same block structure:
\[
\Sigma=\left(\BA{cc}
A & b\\
b^T & c\\
\EA\right),
\]
where $A \in \symm_{(n-1)},~b \in \reals^{(n-1)},~c\in\reals$. In the methods that follow, we will update only one column (and corresponding row) at a time and without loss of generality we can always assume that we are updating the last one.

\paragraph{Block-coordinate descent} The dual problem (\ref{eq:dual-sparseml-block}):
\[
\BA{ll}
\mbox{minimize} & -\log \det (\Sigma + U )-n\\
\mbox{subject to} & |U_{ij}|\leq \rho, \quad i,j=1,\ldots,n.
\EA
\] 
in the variable $U\in\symm_n$, can be written in block format as:
\[
\BA{ll}
\mbox{minimize} & -\log \det (A+V)- \log\left((w+c)-(b+u)^T(A+V)^{-1}(b+u) \right) -n\\
\mbox{subject to} & |w|\leq \rho,~|u_{i}|\leq \rho, \quad i=1,\ldots,n.
\EA
\] 
in the variables $u \in\reals^{(n-1)}$ and $w \in\reals$ ($V$ is fixed at each iteration). We directly get $w=\rho$ so the diagonal of the optimal solution must be $\rho\ones$. The main step at each iteration is then a box constrained quadratic program (QP):
\BEQ\label{eq:box-qp}
\BA{ll}
\mbox{minimize} & (b+u)^T(A+V)^{-1}(b+u)\\
\mbox{subject to} & |u_{i}|\leq \rho, \quad i=1,\ldots,n,
\EA\EEQ
in the variable $u\in\reals^{(n-1)}$. To summarize, the block coordinate descent algorithm proceeds as follows:
\vskip 1ex
\begin{enumerate} \itemsep 0ex
\item Pick the row and column to update.
\item Compute $(A+V)^{-1}$.
\item Solve the box constrained QP in (\ref{eq:box-qp}).
\item Repeat until duality gap is less than precision: 
$\langle \Sigma, X \rangle -n  + \rho \ones^T|X|\ones \leq \epsilon.$
\end{enumerate}
\vskip 1ex
At each iteration, we need to compute the inverse of the submatrix $(A+V)\in\symm_{(n-1)}$, but we can update this inverse using the Sherman-Woodbury-Morrison formula on two rank-two updates, hence it is only necessary to compute a full inverse at the first iteration.

\paragraph{Block-coordinate ascent}

For a fixed $Z$, problem (\ref{eq:sparseml-block}) is equivalent to:
\[\BA{ll}
\mbox{maximize} & \log\left(y-x^TZ^{-1}x\right)-2b^Tx-y(c+\rho)-2 \rho \|x\|_1\\
\mbox{subject to} & y-x^TZ^{-1}x>0,~ y>0,\\
\EA\]
in the variables $x \in \reals^{(n-1)},~y\in\reals$, where $Z \succ 0$ (given) and the Schur complement constraints imply $X \succ 0$. We can solve for the optimal $y$ explicitly and the problem in $x$ becomes:
\[
\max_x -x^TQx-2b^Tx - 2 \rho \|x\|_1,
\]
where $Q:=(c+\rho)Z^{-1}$. Its dual is also box-constrained QP:
\[\BA{ll}
\mbox{minimize} & (b+u)^TZ(b+u)\\
\mbox{subject to} & \|u\|_\infty \leq \rho,
\EA\]
in the variable $u \in \reals^{(n-1)}$. At the optimum for this QP, we must have:
\[
x=-\frac{1}{(c+\rho)}Z(b+u),\quad \mbox{and} \quad y=\frac{1}{(c+\rho)}+\frac{1}{(c+\rho)^2}(b+u)^TZ(b+u).
\]
and we iterate as above.

\paragraph{Smooth optimization for box-constrained QPs}
The two block-coordinate methods detailed in this section both amount to solving a sequence of box-constrained quadratic program of the form:
\BEQ \label{eq:box-qp-section}
\BA{ll}
\mbox{minimize} & x^TAx + b^Tx\\
\mbox{subject to} & \|x\|_\infty \leq \rho,
\EA\EEQ
in the variable $x\in\reals^n$. The objective function has a Lipschitz continuous gradient with constant $L=4c\lambda^{\mathrm{max}}(A)\sqrt{n}$ on the box ${\cal B}=\{x \in \reals^n :~ \|x\|_\infty \leq \rho\}$, where we can define a prox function $(1/2)\|x\|^2$ which is strongly convex with constant 1 and bounded above by $(1/2)n\rho^2$ on ${\cal B}$. From \cite{Nest83} or \cite{Nest03}, we know that solving (\ref{eq:box-qp-section}) up to a precision $\epsilon$ will require at most $2n^{0.75}\sqrt{2\rho^3\lambda^{\mathrm{max}}(A)}/\sqrt{\epsilon}$
iterations of the first-order method detailed in \cite{Nest83}, with each iteration equivalent to a matrix-vector product and a projection on the box ${\cal B}$. This means that the total complexity of solving (\ref{eq:box-qp-section}) is given by:
\[
O\left(n^{2.75}\sqrt{\frac{\rho^3\lambda^{\mathrm{max}}(A)}{\epsilon}}\right)
\]
\paragraph{Complexity estimate} Following \cite{luo92}, with block coordinate descent corresponding to coordinate descent with the almost cyclic rule and using the fact that $\log\det(X)$ satisfies the strict convexity assumptions in \cite[A2]{luo92}, we can show that the convergence rate of the block coordinate descent method is at least linear. Each iteration requires solving a box-constrained QP and takes $O(n^3 log(1/\epsilon))$ operations using an interior point solver or $O(n^{2.75}/\sqrt{\epsilon})$ using the optimal first-order scheme in \cite{Nest83}. We cannot use the same argument to show convergence of block coordinate ascent but empirical performance is comparable. In practice we have found that a small number of sweeps through all columns, independent of problem size $n$, is sufficient for convergence.

\paragraph{Implementation}
The block coordinate descent methods implemented here correspond to coordinate descent using the almost cyclic rule, alternative row/column selection rules could improve the convergence speed. Also, each iteration of the block coordinate descent method corresponds to two rank-two updates of the inverse matrix, hence the cost of maintaining the inverse submatrix using the Sherman Woodbury Morrison formula is only $O(n^2)$.

\section{Numerical results} \label{s:numer}
In this section we test the performance of the methods detailed above on some randomly generated examples. We first form a sparse matrix $A$ with a diagonal equal to one and a few randomly chosen, nonzero off-diagonal terms equal to +1 or -1. We then form the matrix:
\[
B=A^{-1}+\sigma V
\]
where $V\in\symm_n$ is a symmetric, i.i.d uniform random matrix. Finally, we make $B$ positive definite by shifting its eigenvalues, and use this noisy, random matrix to test our covariance selection methods. 

\begin{figure}[ht]
\begin{center}
\psfrag{ap}[c][c]{Noisy inverse $B^{-1}$} \psfrag{bp}[c][c]{Solution for $\rho=0.5$}
\psfrag{cp}[c][c]{Original inverse $A$}
\includegraphics[width= .95 \textwidth]{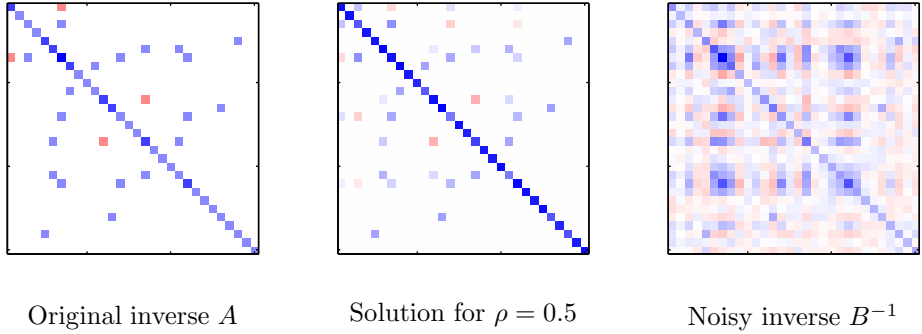}
\caption{\label{fig:pattern} Recovering the sparsity pattern. We plot the original inverse covariance matrix $A$, the solution to problem (\ref{eq:sparseml-relax}) and the noisy inverse $B^{-1}$.}
\end{center}
\end{figure}

\begin{figure}[pht]
\begin{center}
\includegraphics[width= .46 \textwidth]{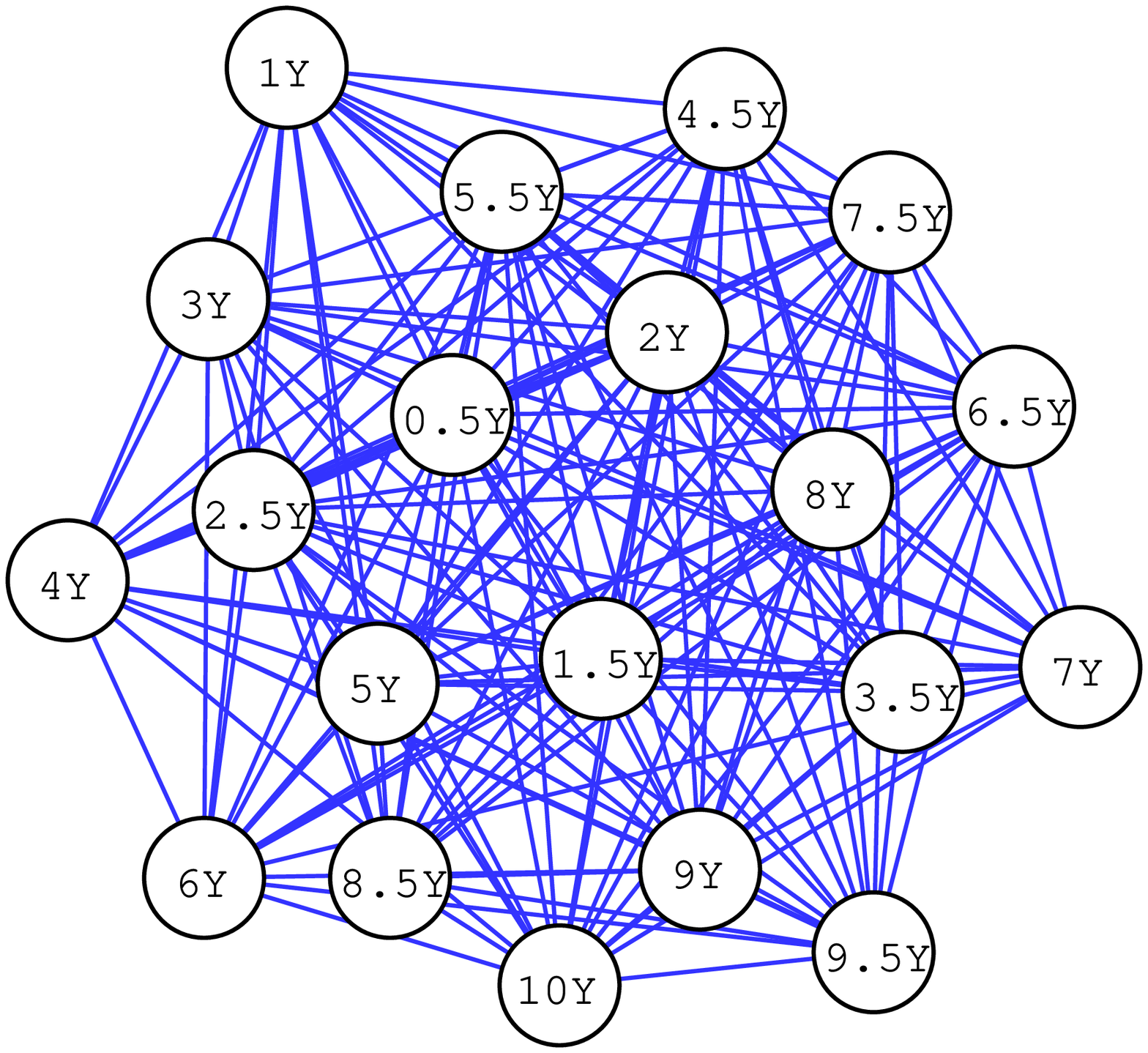}\hfill
\includegraphics[width= .53 \textwidth]{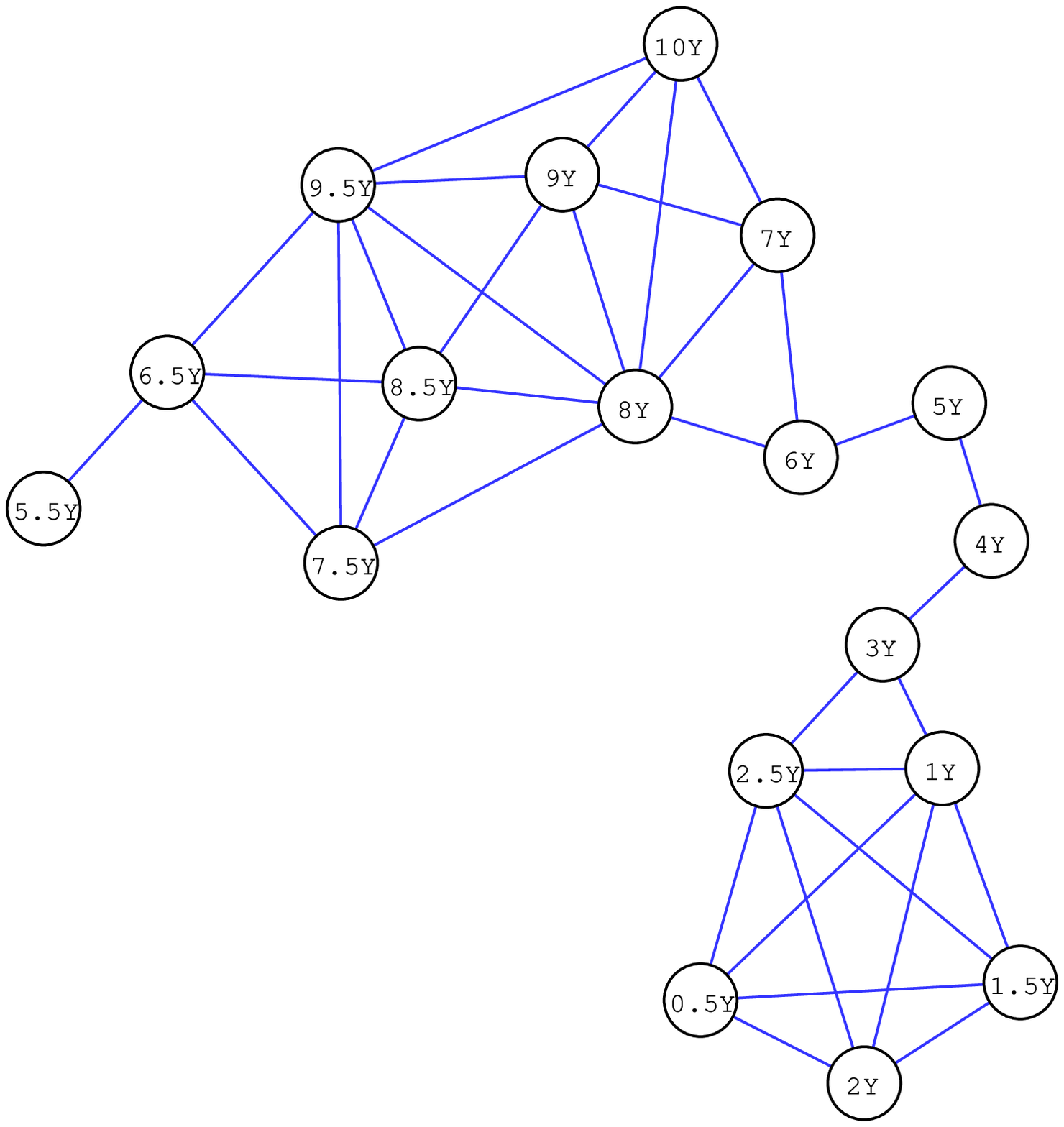}
\caption{\label{fig:rates} We plot the network formed using the solution to problem (\ref{eq:sparseml-relax}) on an interest rate covariance matrix for $\rho=0$ (left) and $\rho=0.1$ (right). In the sparse solution the rates appear clearly clustered by maturity.}
\end{center}
\end{figure}

In Figure \ref{fig:pattern}, we plot the sparsity patterns of the original inverse covariance matrix $A$, the solution to problem (\ref{eq:sparseml-relax}) and the noisy inverse $B^{-1}$ in a randomly generated example with $n=30$, $\sigma=0.15$ and $\rho=0.5$. In Figure \ref{fig:rates} we represent the dependence structure of interest rates (sampled over a year) inferred from the inverse covariance matrix. Each node represents a particular interest rate maturity and the nodes are linked if the corresponding coefficient in the inverse covariance matrix is nonzero (i.e. they are conditionally dependent). We compare the solution to problem (\ref{eq:sparseml-relax}) on this matrix for $\rho=0$ and $\rho=0.1$ and notice that in the sparse solution the rates appear clearly clustered by maturity.

\begin{figure}[pht]
\begin{center}
\psfrag{cpu}[b][t]{CPU Time (in seconds)}
\psfrag{n}[t][b]{Problem Size $n$}
\includegraphics[width=0.48\textwidth]{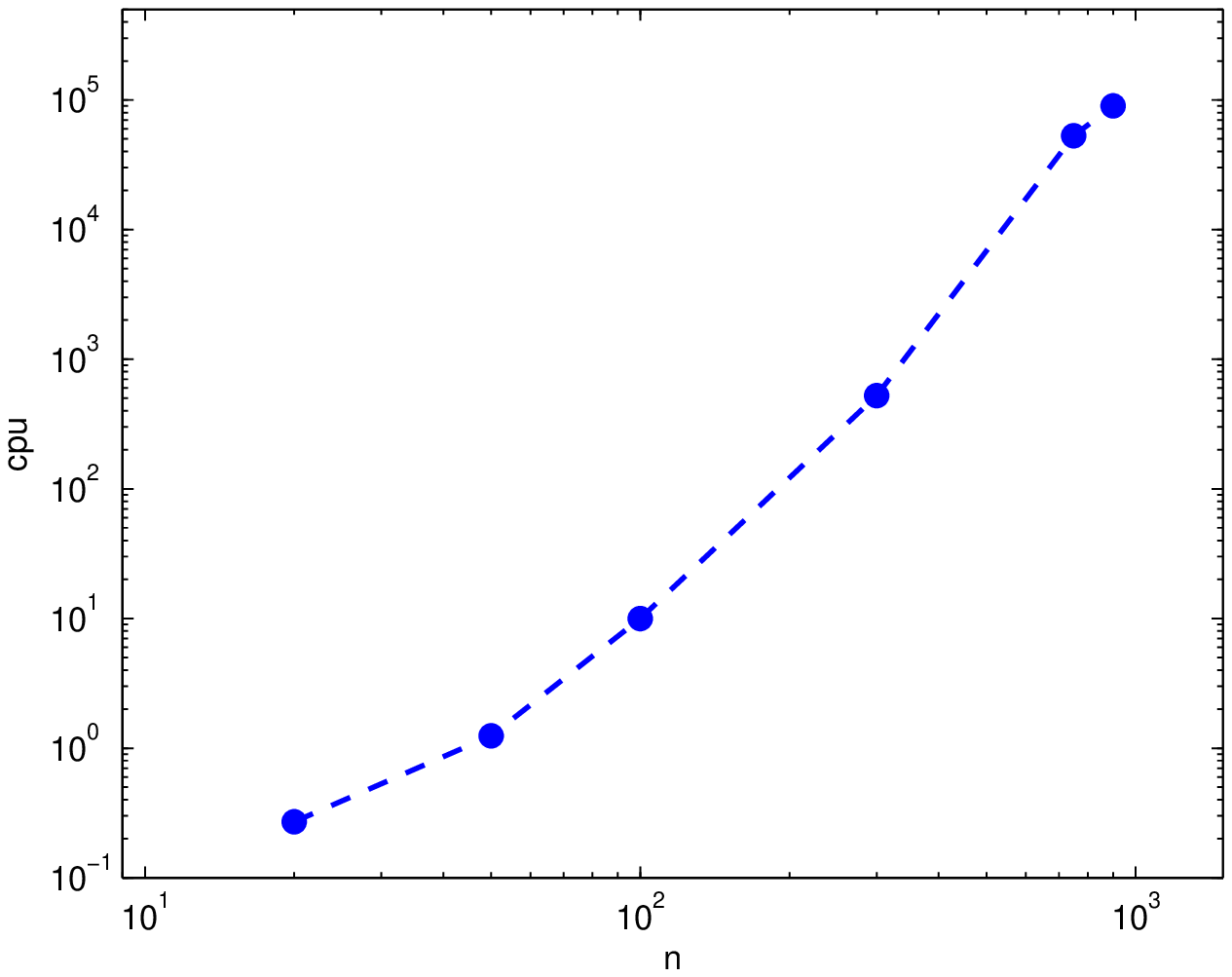}\hfill
\psfrag{cpu}[t][b]{CPU Time (in seconds)}
\psfrag{gap}[b][t]{Duality Gap}
\includegraphics[width=0.48 \textwidth]{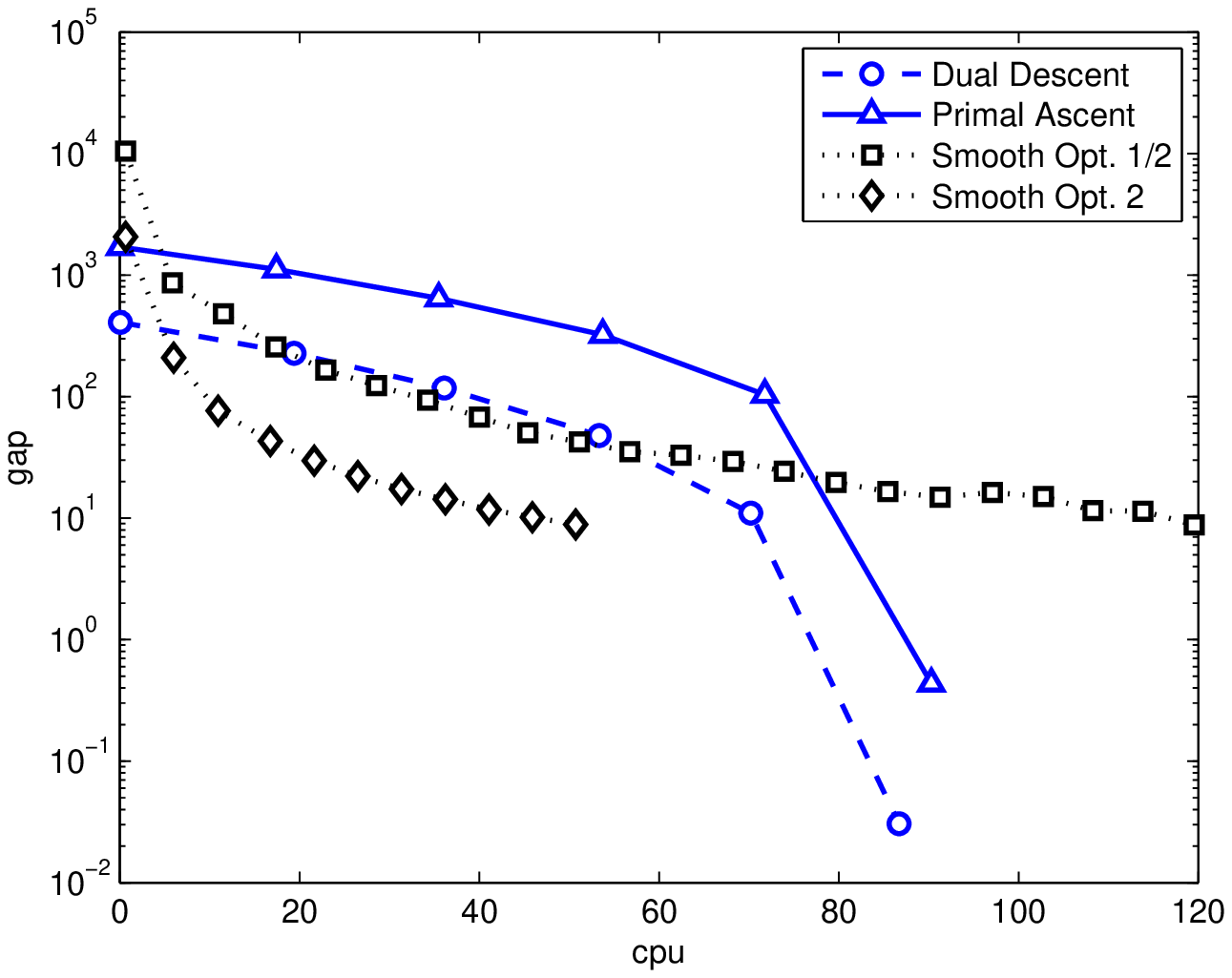}
\caption{\label{fig:cputime} Computing time. \emph{Left}: We plot CPU time to reduce the duality gap by a factor $10^{-2}$ versus problem size $n$, on randomly generated problems, using the coordinate descent code and the optimal first-order algorithm for solving box QPs. \emph{Right}: We plot duality gap versus CPU time for both smooth minimization and block-coordinate algorithms, for a problem of size $n=250$.}
\end{center}
\end{figure}

In Figure \ref{fig:cputime}, we study computing times for various choices of algorithms and problem sizes. On the \emph{left}, we plot CPU time to reduce the duality gap by a factor $10^{-2}$ versus problem size $n$, on randomly generated problems, using the coordinate descent code and the optimal first-order for solving box QPs. On the \emph{right}, we plot duality gap versus CPU time for both smooth minimization and block-coordinate algorithms for a randomly generated problem of size $n=250$. For the smooth minimization code, we set $\alpha=1/\lambda^{\mathrm{max}}(B)$ and we plot computing time for both $\beta=1/(2\lambda^{\mathrm{min}}(B))$ (Smooth. Opt. 1/2) and $\beta=2/\lambda^{\mathrm{min}}(B)$ (Smooth. Opt. 2).

\section*{Acknowledgments.} The authors would like to thank Francis Bach, Peter Bartlett and Martin Wainwright for enlightening discussions on the topic. We would also like to acknowledge financial support from NSF grant DMS-0625352, EUROCONTROL grant C20083E/BM/05 and a gift from Google, Inc.

\small{
\bibliographystyle{siam}
\bibliography{MainPerso}}
\end{document}